\documentclass{amsart}

\newcommand{\RR}{\mathbb{R}} 
\newcommand{\CC}{\mathbb{C}} 

\newcommand{\HP}{\mathbf{H}}

\newcommand{\RE}{\mathrm{Re}\,}

\newcommand{\Aut}{\mathrm{Aut}}

\newcommand{\dil}{\mathcal{D}} 
\newcommand{\tran}{\mathcal{T}} 

\newcommand{\JB}{\mathcal{N}}

\newcommand{\UU}{\mathrm{U}}

\newcommand{\Dil}{\mathcal{D}}
\newcommand{\Hei}{\mathcal{H}}

\newcommand{\rot}{\mathcal{R}}

\newcommand{\cC}{\mathcal{C}}

\newcommand{\paren}[1]{\left(#1\right)}
\newcommand{\abs}[1]{\left\lvert#1\right\rvert}
\newcommand{\norm}[1]{\left\|#1\right\|}

\newcommand{\ip}[1]{\left\langle#1\right\rangle}
\newcommand{\pd}[2]{\frac{\partial#1}{\partial#2}}

\newtheorem{theorem}{Theorem}[section]

\newtheorem{lemma}[theorem]{Lemma}

\theoremstyle{definition}
\newtheorem{definition}[theorem]{Definition}
\theoremstyle{remark}
\newtheorem{remark}[theorem]{Remark}

\numberwithin{equation}{section}

\title[Homogeneous almost complex manifolds \& 
compact quotients] {Homogeneous almost complex manifolds and 
their compact quotients}

\author{Kang-Tae Kim, Kang-Hyurk Lee and Yoshikazu Nagata}

\address {Kang-Tae Kim and Yoshikazu Nagata: Department of 
Mathematics and Center for Geometry and its Applications, 
Pohang University of Science and Technology, 37673, The 
Republic of Korea} 
\email{kimkt@postech.ac.kr, yoshikazu@postech.ac.kr}
\address{Kang-Hyurk Lee: Department of Mathematics and 
Research Institute of Natural Science, Gyeongsang National 
University, Jinju, Gyeongnam, 660-701, The Republic of Korea}
\email{nyawoo@gnu.ac.kr}

\thanks{Research of the authors is supported in part by the 
NRF Grant 2011-0030044 (SRC-GAIA) of The Republic of Korea.}

\subjclass[2010]{32M05}

\keywords{$J$-automorphism, compact quotient}

\begin{document}

\maketitle

\begin{abstract}
This paper investigates the (non)existence of compact quotients of 
the homogeneous almost-complex strongly-pseudoconvex manifolds 
discovered and classified by Gaussier-Sukhov 
\cite{Gau-Sukh, Gau-Sukh2} and K.-H. Lee \cite{KHL1, KHL2}.
\end{abstract}

\section{Introduction}

Let \(M\) be an almost complex manifold of real dimension \(2m\), 
\(m\ge 1\), with an almost complex structure \(J\).  It is said to
be modeled after bounded strictly pseudoconvex 
domains, if the following two properties hold:
\begin{itemize}
\item[(1)] \(M\) is $J$-holomorphically equivalent to a subdomain of 
another almost complex manifold with its Levi form (cf. Section 2) at
every boundary point positive-definite, and 
\item[(2)] \(M\) is Kobayashi hyperbolic.
\end{itemize}
It was generally believed for some time that, if one further 
assumes for it to be \textit{homogeneous}, meaning that the action by  
the group \(\Aut(M, J)\) of all \(J\)-holomorphic
diffeomorphisms of \(M\) into itself is transitive, then 
the almost complex structure should be integrable, and 
consequently---due to the well-known theorem by Wong and 
Rosay---the complex manifold \((M,J)\) would have to be 
biholomorphic to the unit ball \(B^m\) in \(\CC^m\).

While such belief was justified in the case of \(m\le 2\) 
by Gaussier and Sukhov \cite{Gau-Sukh, Gau-Sukh2}, they 
in contrast showed that there exists an example 
indicating that it is not the case if \(m>2\).  Then the second named 
author of this article classified all such manifolds for every
\(m \ge 3\) in \cite{KHL1, KHL2} (cf.\ Theorem \ref{KHL}). 
It has turned out that there are infinitely many such 
examples, $J$-holomorphically inequivalent to each other.
\smallskip

Upon such observation, there arises a natural question whether 
the homogeneous manifolds (noncompact) obtained here would 
admit a compact quotient by a discrete subgroup.  This question was 
asked to the first named author several different times by many 
prominent mathematicians, from at least 10 years ago, 
when he gave lectures explaining this line of research at the 
institutions including l'\'Ecole Polytechnique de Palaiseau of France, 
Peking University of China, the Korea Institute 
for Advanced Study (Seoul) of Korea, and other places. It is our
pleasure to acknowledge our indebtedness.
\medskip

The purpose of this article is, indeed, to provide the answer---negative,
however.  Deferring the introduction of terminology and necessary 
definitions to the later sections, we present our results first.

\begin{theorem} \label{main}
Let \(\Omega\) be an open connected subset 
of an almost complex manifold \((M, J)\) of real dimension \(2m\)
with \(\cC^2\) smooth strictly pseudoconvex boundary.
If there exists a discrete subgroup \(\Gamma\) of the automorphism
group  \(\Aut(\Omega, J)\) such that \(\Omega /\Gamma\) is 
compact, then \((\Omega, J)\) is biholomorphic to the standard
unit ball in \(\CC^m\).
\end{theorem}

From here on, open connected subsets of a manifold will be
called {\it domains}, as usual. Now, notice that Theorem 
\ref{main} implies, according to the discussion above, the following 
result concerning the nonexistence of compact quotients.

\begin{theorem} \label{main-coro}
If \(\Omega\) is a homogeneous domain
with \(\cC^2\) smooth strictly pseudoconvex boundary 
in an almost complex manifold \((M, J)\) of real dimension \(2m\) 
whose \(J\)-structure is non integrable, then \((\Omega, J)\) does
not admit any compact quotient by any discrete subgroup of
\(\Aut (\Omega, J)\). 
\end{theorem}

We remark that these theorems are significant especially for 
the case \(m \ge 3\).

\section{Model manifolds/domains}

For a domain \(\Omega\) in \( (M,J)\), any \(\mathcal C^\infty\) smooth diffeomorphism 
\(\psi \colon \Omega \to \Omega\) is called a \(J\)-\textit{holomorphic 
automorphism} if \(J \circ d\psi = d\psi \circ J\).  (More generally, 
a smooth map \(f\colon (M,J) \to (\tilde M, \tilde J)\) between two almost 
complex manifolds is called \textit{pseudo-holomorphic} or, more precisely, 
\((J,\tilde J)\)-\textit{holomorphic} if \(\tilde J \circ df = df \circ J\) holds.)
Such automorphisms
form a topological group, denoted by \(\Aut (\Omega, J)\)
under the law of composition endowed with
the compact-open topology.

Moreover, if the boundary $\partial\Omega$ of $\Omega$ is 
\(\mathcal C^2\) smooth, then the implicit function theorem implies 
that for every boundary point $p$ there are an open neighborhood 
\(V\) of \(M\) and a \(\mathcal C^2\) function \(\rho\colon V \to \RR\) 
such that \(\Omega \cap V = \{ z \in V \colon \rho(z)<0\} \) and that 
\(d\rho(q) \neq 0\) for any \(q \in \partial\Omega \cap V\).  Such 
\(\rho\) is called a \textit{local defining function}.  

For a 1-form \(\omega\), the dual \(J^*\) of \(J\) is defined by 
\((J^* \omega)(v) = \omega (Jv)\).  Then the \textit{Levi form}
of \(\rho\) is defined to be 
\[
\mathcal L_\rho (v, w) := -d(J^* d\rho) (v, Jw).
\]
Then we say that \(\Omega\) is \textit{strictly $J$-pseudoconvex} at 
\(p \in \partial\Omega\) if 
\( \mathcal L_\rho (v, v) > 0 \) for every nonzero vector 
\(v \in T_p\partial\Omega \cap JT_p\partial\Omega\).  
\medskip

Now write \(m=n+1 \ge 1\), let $z^0,z^1,\ldots,z^n$ 
represent the standard coordinate functions of $\CC^m = \CC^{n+1}$, 
and denote by $z'=(z^1,\ldots,z^n)$ the standard coordinate system 
of $\CC^n$.
\smallskip

For a mapping $\phi$ to $\CC^{n+1}$, denote by 
$\phi^0=z^0\circ\phi$ and $\phi'=z'\circ\phi$, and hence 
$\phi=(\phi^0,\phi')$. Greek indices $\alpha,\beta,\ldots$ run from 
$1$ to $n$ and the summation convention is always assumed: 
$\JB_{\alpha\beta}z^\beta=\sum_{\beta=1}^n\JB_{\alpha\beta}z^\beta$,
for instance.
We also put bar on the indices to denote the complex conjugation of 
the corresponding tensor coefficients such as: 
$\JB_{\bar\alpha\bar\beta}=\overline{{\JB}_{\alpha\beta}}$, 
$A^{\bar\alpha}_{\bar\beta}=\overline{A^\alpha_\beta}$. 

\medskip

Let $\HP=\{(z^0,z')\in\CC\times\CC^n: \RE z^0 +\norm{z'}^2<0\}$ 
be the Siegel half space. But we shall endow an almost complex structure
potentially different from the standard integrable one. 
\smallskip

For each $n\times n$ skew-symmetric matrix $\JB=(\JB_{\alpha\beta})$, 
define the almost complex structure $J_\JB$ of $\CC^{n+1}$ by
\begin{align*}
	 J_\JB &= i \pd{}{z^0} \otimes d z^0 +
	 	\paren{i\pd{}{z^\alpha} 
		+
		2\JB_{\alpha\beta}z^\beta\pd{}{z^{\bar0}}} \otimes d z^\alpha \\
		&\quad
	-i\pd{}{z^{\bar0}}  \otimes d z^{\bar0} 
		+ \paren{ -i\pd{}{z^{\bar\alpha}} 
		+2\JB_{\bar\alpha\bar\beta}z^{\bar\beta}\pd{}{z^0} } 
		\otimes d z^{\bar\alpha} \;.
\end{align*}
Note that the almost complex structure $J_\JB$ is integrable if and 
only if $\JB=0$.
\medskip

Then the following characterization of the strictly $J$-pseudoconvex
domains with an automorphism orbit accumulating at a boundary point
has been established earlier, by the second named author of this article:

\begin{theorem}[\cite{KHL2}] \label{KHL}
If \((\Omega, J)\) is a domain in an almost complex manifold \((M,J)\) with  
\(\mathcal{C}^2\) smooth strictly pseudoconvex boundary point 
$q \in \partial\Omega$ such that there are a point $p \in \Omega$ and 
a sequence of $J$-automorphisms \(\varphi_j \in \Aut (\Omega, J)\) 
satisfying \(\displaystyle \lim_{j\to\infty} \varphi_j(p)=q\), 
then \( (\Omega, J) \) is \((J, J_\JB)\)-biholomorphic to one of 
\( (\HP,J_\JB) \).  
\end{theorem}

\begin{remark} \rm
Every \( (\HP,J_\JB) \) is homogeneous, as one sees in the next
section.
\end{remark}

\section{Automorphisms of \( (\HP,J_\JB) \)} 

Note that, due to the preceding discussion, the proof of 
Theorem~\ref{main} reduces to demonstrating the (non)existence 
of the compact quotients for the domains \( (\HP,J_\JB) \), 
called the \textit{model domains} in \cite{KHL1, KHL2}.
\medskip

There are four types of automorphisms generating the automorphism
group of the model domain $(\HP,J_\JB)$ for $\JB\neq 0$:

\begin{itemize}
\item[(1)] \textbf{$\RR$-action:} For each $s\in\RR$,
\begin{equation*}
\tran_s
	=(z^0+is,z')\in\Aut(\HP,J_\JB)
	\;.
\end{equation*}

\item[(2)] \textbf{$\CC^n$-action:} For each $w'\in\CC^n$, let us define
\begin{equation*}
h_{w'}(z')=
	-\norm{w'}^2-2\ip{z',\bar w'}
	+i(\JB_{\alpha\beta}z^\alpha w^\beta 
	+\JB_{\bar\alpha\bar\beta}z^{\bar\alpha}w^{\bar\beta})
	\;.
\end{equation*}
Here $\ip{\,\cdot\,,\,\cdot\,}$ is the standard hermitian product of 
$\CC^n$: $\ip{z',\bar w'}=\delta_{\alpha\bar\beta} z^\alpha 
w^{\bar\beta}$. Then $(\CC^n,+)$ can be embedded into 
$\Aut(\HP,J_\JB)$ by
\begin{equation*}
w'\quad\longmapsto\quad\Hei_{w'}
	=(z^0+h_{w'}(z'), z'+w')\in\Aut(\HP,J_\JB) 
	\;,
\end{equation*}
since $\Hei_{w'}\circ\Hei_{v'}=\Hei_{w'+v'}$ for any $w',v'\in\CC^n$.
\end{itemize}
\medskip

We remark in passing that the subgroup generated by $\Hei_{w'}$ and 
$\tran_s$ is in fact isomorphic to the Heisenberg group.
\medskip

\begin{itemize}
\item[(3)] \textbf{Dilation:} For each $t>0$, we have
\begin{equation*}
\dil_t=(tz^0,t^{1/2}z')\in\Aut(\HP,J_\JB) \;.
\end{equation*}

\item[(4)] \textbf{Isotropy:} Let $\UU_\JB$ be the set of $n\times n$ 
complex matrices $A=(A^\alpha_\beta)$ with
\begin{equation*}
\delta_{\alpha\bar\beta} 
	=A^\mu_\alpha \delta_{\mu\bar\nu}A^{\bar\nu}_{\bar\beta}
	\quad\text{and}
\quad \JB_{\alpha\beta} 
	=A^\mu_\alpha\JB_{\mu\nu}A^\nu_\beta
	\;,
\end{equation*}
i.e., $I=A^t \overline{A}$ and $\JB= A^t \JB A$. Then $\UU_\JB$ can be
realized as a compact subgroup of the unitary group $\UU(n)$ which 
can be embedded into $\Aut(\HP,J_\JB)$ via
\begin{equation*}
A\quad\longmapsto\quad \rot_A=(z^0,Az')\in\Aut(\HP,J_\JB) \;.
\end{equation*}
\end{itemize}
\bigskip

Then we have:

\begin{theorem}[\cite{KHL2}]
For each $\phi\in\Aut(\HP,J_\JB)$ with \(\JB \neq 0\), there is 
a unique choice for 
$t > 0$, $s\in\RR$, $w'\in\CC^n$ and $A\in\UU_\JB$ such that
\begin{equation*}
\phi=\tran_s\circ\Hei_{w'}\circ\dil_t\circ\rot_A
=\paren{
	tz^0+ h_{w'}(z')+is, 
	t^{1/2}Az'+w'
	}
	.
\end{equation*}
\end{theorem}

\section{Discrete subgroups and the limit sets}

In this section, we consider only the case \(\JB \neq 0\).

\begin{definition}
Let \(H\) be a subgroup of the automorphism group \(\Aut(\HP,J_\JB)\). 
By the \textit{limit set} of \(H\) we mean the set \(\Lambda(H)\) of 
all accumulation points of the orbits by \(H\).  Here, the limit set may 
contain the points at infinity.
\end{definition}

Suppose that $\Gamma$ is a discrete subgroup of $\Aut(\HP,J_\JB)$. 
The aim of this section is to analyze the limit set 
\(\Lambda(\Gamma)\), which will eventually lead us to 
the proof of Theorem \ref{main}.
\medskip

Choose $\phi_\nu\in\Gamma$ and write
\begin{equation*}
\phi_\nu (z)=
	\paren{
	t_\nu z^0+ h_{w'_\nu}(z')+is_\nu, 
	t_\nu^{1/2}A_\nu z'+w'_\nu
	}
\end{equation*}
for some $t_\nu >0$, $s_\nu\in\RR$, $w'_\nu\in\CC^n$ and 
$A_\nu\in\UU_\JB$. 

For $p=(p^0,p')\in\HP$, we assume that $\phi_\nu(p)\to(q^0,0')$ as 
$\nu\to\infty$, where $q^0\in\{is:s\in\RR\}\cup\{\infty\}$. Since 
$A_\nu$ is an element of the compact group $\UU_\JB$, we may 
assume that $A_\nu\to A_\infty$ in $\UU_\JB$. We may assume 
that $t_\nu$ converges in $[0,+\infty]$, and also $s_\nu$ in
$[-\infty, +\infty]$, respectively.
Note that $w'_\nu$ can also be assumed to converge in the 
one-point compactification of $\CC^n$.
\bigskip

\noindent\textbf{Case (1):} $\norm{w_\nu'}\to\infty$. Since 
$t_\nu^{1/2}A_\nu p'+w'_\nu\to 0'\in\CC^n$ and $A_\nu\to A_\infty$, 
we have $t_\nu\to\infty$ and $t_\nu^{-1/2}w'_\nu\to - A_\infty p'$. 
Let us consider
\begin{multline*}
\RE\phi_\nu^0(p)=t_\nu\RE p^0 -\norm{w'_\nu}^2
-2\RE\ip{p',\bar w'_\nu}
	\\
	=t_\nu\paren{
	\RE p^0 -\norm{t_\nu^{-1/2}w'_\nu}^2-2\RE\ip{t_\nu^{-1/2}p',
	t_\nu^{-1/2} 	\bar w'_\nu}
	}
	.
\end{multline*}
Since $\RE p^0 -\norm{t_\nu^{-1/2}w'_\nu}^2-2\RE\ip{t_\nu^{-1/2}p',
t_\nu^{-1/2}\bar w'_\nu}\to \RE p^0 -\norm{p'}^2<0$, we have 
$\phi_\nu^0(p)\to\infty$. Thus the accumulating point is $(\infty,0')$.

\bigskip

\noindent\textbf{Case (2):} $w'_\nu\to w'_\infty\neq 0'$. Since 
$t_\nu^{1/2}A_\nu p'+w'_\nu\to 0'\in\CC^n$ and $A_\nu\to A_\infty$, 
we have $t_\nu\to t>0$. If $s_\nu\to s\in\RR$, then
\begin{multline*}
\phi_\nu(z)=
	\paren{
	t_\nu z^0+ h_{w'_\nu}(z')+is_\nu, 
	t_\nu^{1/2}A_\nu z'+w'_\nu
	}
	\\
	\longrightarrow\quad
	\phi(z)=\paren{
	tz^0+ h_{w'_\infty}(z')+is, 
	t^{1/2}A_\infty z'+w'_\infty
	}\in\Aut(\HP,J_\JB)
	\;.
\end{multline*}
This is a contradiction to the discreteness of $\Gamma$. Hence 
$\abs{s_\nu}\to\infty$ so that the accumulation point is $(\infty,0')$ only.

\bigskip

\noindent\textbf{Case (3):} $w'_\nu\to 0'$. Suppose that $t_\nu$ does not 
converge to $0$. If $p'\neq 0'$, then $t_\nu\to 0$ because 
$t_\nu^{1/2}A_\nu p'+w'_\nu\to 0'\in\CC^n$. Thus we must have 
$p'=0'$ and $\RE p^0<0$. 

If $t_\nu\to\infty$, then
\begin{equation*}
\RE\phi_\nu^0(p)=t_\nu\RE p^0 -\norm{w'_\nu}^2-2\RE\ip{p',
\bar w'_\nu} = t_\nu\RE p^0 -\norm{w'_\nu}^2
\end{equation*}
and this implies that $\RE\phi_\nu^0(p)\to-\infty$. 
So the accumulation point is $(\infty,0')$.

If $t_\nu\to t>0$, by the same argument as in Case (2) 
$\abs{s_\nu}\to\infty$. 
Thus the accumulation point is $(\infty,0')$.

\medskip

It remains to analyze the case of $t_\nu\to 0$. For this purpose, 
we pose the following

\begin{lemma}\label{lem:small t}
Let $\Gamma$ be a discrete subgroup of $\Aut(\HP,J_\JB)$. If there 
is $\phi\in\Gamma$ of the form
\( \displaystyle 
\phi(z)=
	\paren{tz^0+is,t^{1/2}Az'}
\)
for some $t<1$, then each $\psi\in\Gamma$ is of the form
\( \displaystyle
\psi(z)=
	\paren{az^0+ib,a^{1/2}Bz'}
\)
with
\( \displaystyle
\frac{b}{1-a}=\frac{s}{1-t}
\).
\end{lemma}
\smallskip

Postponing the proof of this lemma to the end of this section, we 
continue analyzing the accumulation points.
\bigskip

Assume that $t_\nu\to 0$ as $\nu \to\infty$. Then fix 
$\phi_m$ such that $t_m<1$ and, consequently, 
$I-\sqrt{t_m}A_m$ is invertible. From now on, we simply denote 
by $t=t_m$, $s=s_m$, $w'=w'_m$, 
$A=A_m$ and put $\zeta'=(I-\sqrt{t}A)^{-1}w'$ and 
$\gamma=\Hei_{\zeta'}\circ\Dil_t\circ\rot_A$.  Thus we have
\begin{align*}
\gamma'(z')&=t^{1/2}Az'+\zeta' \\
(\gamma^{-1})'(z') &= t^{-1/2}A^{-1}(z'-\zeta'),
\end{align*}
which immediately implies
\begin{equation}\label{eqn:inner}
\gamma^{-1}\circ\phi_m\circ\gamma(z)=
	\paren{
	tz^0+i\tilde s,
	t^{1/2}Az'
	},
\end{equation}
for some $\tilde s\in\RR$. For general integer values of $\nu$, 
it follows that
\begin{multline*}
(\gamma^{-1}\circ\phi_\nu\circ\gamma)'(z')
	=
	(\gamma^{-1})'\circ\phi_\nu'\circ\gamma'(z')
	\\
	=
	t_\nu^{1/2}A^{-1}A_\nu Az'
	+t^{-1/2}t_\nu^{1/2}A^{-1}A_\nu\zeta'
	+t^{-1/2}A^{-1}w'_\nu
	-t^{-1/2}A^{-1}\zeta'.
\end{multline*}
Since $\gamma^{-1}\circ\phi_m\circ\gamma$ is of the same form 
as $\phi$ in Lemma~\ref{lem:small t}, we obtain, as far as the element 
$\gamma^{-1}\circ\phi_\nu\circ\gamma$ of the discrete group 
$\gamma^{-1}\Gamma\gamma$ is concerned, that
\begin{equation*}
(\gamma^{-1}\circ\phi_\nu\circ\gamma)'(z') = t_\nu^{1/2}C_\nu z',
\end{equation*}
where $C_\nu=A^{-1}A_\nu A$.  Now we have
\begin{align*}
(\phi_\nu)'(z') 
	& = 
	\gamma'\circ(\gamma^{-1}\circ\phi_\nu\circ
	\gamma)'\circ(\gamma^{-1})'(z)
	\\
	& =
	t_\nu^{1/2}AC_\nu A^{-1}z'-t_\nu^{1/2}AC_\nu A^{-1}\zeta'+\zeta' 
	~~\longrightarrow~~ \zeta'
\end{align*}
as \(\nu \to \infty\) for any $z'\in\CC^n$. Since $(\phi_\nu)'(p') \to 0'$ 
as $\nu\to\infty$, it follows that $\zeta'=0$, which in turn implies
$w'=w_m'=0$. Thus $\phi_m(z)=(tz^0+is, t^{1/2}Az')$ with $t<1$. 
Applying 
Lemma~\ref{lem:small t} again for $\phi=\phi_m$, each 
$\psi\in\Gamma$ is of the form
\begin{equation*}
\psi(z)=
	\paren{az^0+ib,a^{1/2}Bz'}
\end{equation*}
with
\begin{equation*}
\frac{b}{1-a}=\frac{s}{1-t}\;.
\end{equation*}
Especially,
\begin{equation*}
\phi_\nu (z)=
	\paren{
	t_\nu z^0+is_\nu, 
	t_\nu^{1/2}A_\nu z'
	}
\end{equation*}
where
\begin{equation*}
\frac{s_\nu}{1-t_\nu}=\frac{s}{1-t}\;.
\end{equation*}
for any $\nu$. Thus
\begin{equation*}
s_\nu=\frac{s(1-t_\nu)}{1-t}\longrightarrow\frac{s}{1-t}
\end{equation*}
so
\begin{equation*}
\phi_\nu(z)=
	\paren{
	t_\nu z^0+is_\nu, 
	t_\nu^{1/2}A_\nu z'
	}
	\to
	(is/(1-t),0').
\end{equation*}
Since $s/(1-t)$ is an invariant of $\Gamma$, the accumulation point 
should be $(is/(1-t),0')$ when $w'_\nu\to 0'$ and $t_\nu\to0$.

\bigskip

\noindent\textit{Proof of Lemma~\ref{lem:small t}.}
Since $\phi^{-1}(z) = \paren{t^{-1}z^0 -t^{-1}is, t^{-1/2}A^{-1}z'}$, we have
\begin{align*}
\phi^k(z) 
	&= 
	\paren{
	t^kz^0 +  is\sigma_k, 
	t^{k/2}A^kz'
	} \\
\phi^{-k}(z) 
	&= 
	\paren{
	t^{-k}z^0 -ist^{-k}\sigma_k, 
	t^{-k/2}A^{-k}z'
	}
\end{align*}
where $\phi^k$ is the $k$-th iteration of $\phi$  and $\phi^{-k}=(\phi^k)^{-1}$ and 
\begin{equation*}
\sigma_k=1+t+\cdots+t^{k-1}\;.
\end{equation*} 
Take $\psi\in\Gamma$ and let
\begin{equation*}
\psi(z)=
	\paren{
	az^0+h_{v'}(z') +ib, 
	a^{1/2}Bz'+v'
	}
\end{equation*}
Then
\begin{multline*}
\lefteqn{\phi^k\circ\psi\circ\phi^{-k}(z)} \\
	=
	\paren{
	az^0 +i\paren{(1-a)s\sigma_k+t^kb}
	+t^kh_{v'}(t^{-k/2}A^{-k}z') , 
	a^{1/2}A^kBA^{-k}z'+t^{k/2}A^kv'
	}\;.
\end{multline*}
Since $t^kh_{v'}(t^{-k/2}A^{-k}z')\to 0 $ as $k\to\infty$, 
$\phi^k\circ\psi\circ\phi^{-k}$ subsequentially converges to 
some $\tilde\psi\in\Gamma$ of the form
\begin{equation*}
\tilde\psi= \paren{
	az^0 +(1-a)is\sigma  , 
	a^{1/2}\widetilde{A}B\widetilde{A}^{-1}z'
	}
\end{equation*}
where $\sigma=\lim \sigma_k=1/(1-t)$. Since $\Gamma$ is discrete, 
we conclude that $v'=0$ and
\begin{equation*}
(1-a)s\sigma_k+t^kb = (1-a)s\sigma
\end{equation*}
for infinitely many $k$. Thus
\begin{equation*}
\frac{b}{1-a}=\frac{1}{t^k}s(\sigma-\sigma_k) = s\sigma=\frac{s}{1-t}\;.
\end{equation*}
This completes the proof. \qed

\section{Proof of Theorems \ref{main} and \ref{main-coro}}

Let \(\Gamma\) be a discrete subgroup of  \(\Aut(\Omega, J)\) with
\(\Omega/\Gamma\) compact. Then there exists a (finite) boundary point
at which the orbit of \(\Gamma \) accumulates. By Theorem
\ref{KHL}, \(\Omega\) is biholomorphic to a model \( (\HP,J_\JB) \), say, 
that is homogeneous and complete Kobayashi hyperbolic.
Suppose now that \(\Omega\), or equivalently \( (\HP,J_\JB) \),  
is not biholomorphic to the standard unit ball of \(\CC^m\).

According to the arguments of Section 4 above, the set 
\(\Lambda(\Gamma) \cap (\CC\times\{(0,\ldots,0)\})\) consists of  
at most two elements: the \textit{extended boundary point}
\((\infty, 0, \ldots, 0)\), or the boundary point of type 
\((a,0,\ldots, 0)\) for some \(a \in \CC\).  
Suppose that \( (\HP,J_\JB) \) admits a compact quotient
by a discrete group \(\Gamma\), say. Then, since \( (\HP,J_\JB) \) is 
complete Kobayashi hyperbolic, we must have  
\[
\Lambda(\Gamma) \cap (\CC\times\{(0,\ldots,0)\}) 
= \partial\HP\cap (\CC\times\{(0,\ldots,0)\}) =
\RR\times\{(0,\ldots,0)\},
\] 
which cannot coincide with any two element set.
This yields the proof of Theorems \ref{main} as well as \ref{main-coro}.
\hfill
$\Box$


\begin{thebibliography}{99}

\bibitem{Gau-Sukh} H. Gaussier; A. Sukhov: Wong–Rosay theorem in 
almost complex manifolds, \texttt{http:\\www.arXiv.org: 
math.CV/0307335}.

\bibitem{Gau-Sukh2} H. Gaussier; A. Sukhov: On the geometry of 
model almost complex manifolds with boundary, \textit{Math. Z.} 254
(2006), no. 3, 567--589.

\bibitem{KHL1} K.-H. Lee: Domains in almost complex manifolds with 
an automorphism orbit accumulating at a strongly pseudoconvex 
boundary point. \textit{Michigan Math. J.} 54 (2006), no. 1, 179--205. 

\bibitem{KHL2} K.-H. Lee:  Strongly pseudoconvex homogeneous 
domains in almost complex manifolds. \textit{J. Reine Angew. Math.} 
623 (2008), 123--160.

\bibitem{Ros} J.-P. Rosay: Sur une caractérisation de la boule parmi 
les domaines de ${\bf C}^{n}$ par son groupe d'automorphismes. (French) 
\textit{Ann. Inst. Fourier (Grenoble)} 29 (1979), no. 4, ix, 91--97.

\bibitem{Wong} B. Wong: Characterization of the unit ball in 
$\CC^n$ by its automorphism group. \textit{Invent. Math.} 41 (1977), 
no. 3, 253--257.
\end{thebibliography}
\end{document}